\documentclass{amsart}
\usepackage{amssymb,amscd}

\usepackage[french]{babel}
\usepackage{times}

\usepackage[latin1]{inputenc}
\usepackage[T1]{fontenc}
\newcounter{intro}

\newtheorem{theo}[intro]{Théorème}

\newtheorem{thm}{Théorème}[section]
\newtheorem{lem}[thm]{Lemme}
\newtheorem{prop}[thm]{Proposition}

\theoremstyle{remark}
\newtheorem{rem}[thm]{Remarques}

\numberwithin{equation}{section}

\newcommand{\cref}[1]{Corollary~\ref{#1}}

\newcommand{\C}{\mathbb{C}}

\newcommand{\Ag}{\mathfrak{a}}
\newcommand{\Gg}{\mathfrak{g}}

\newcommand{\Pg}{\mathfrak{p}}

\newcommand{\Ng}{\mathfrak{n}}

\newcommand{\Kg}{\mathfrak{k}}

\renewcommand{\Re}{\rm Re\,}

\DeclareMathOperator{\ad}{ad}

\DeclareMathOperator{\Id}{Id}

\DeclareMathOperator{\vol}{vol}

\newcounter{counteroman}
\newenvironment{enumeroman}{\begin{list}{\roman{counteroman})}{\usecounter{counteroman}}}{\end{list}}

\begin{document}

\title[Estimées des noyaux de Green et de la chaleur sur les espaces symétriques. ]
{Estimées des noyaux de Green et de la chaleur sur les espaces symétriques.}

\author{Gilles Carron}
\address{Laboratoire de Math\'ematiques Jean Leray (UMR 6629), Universit\'e de Nantes, 
2, rue de la Houssini\`ere, B.P.~92208, 44322 Nantes Cedex~3, France}
\email{Gilles.Carron@math.univ-nantes.fr}


\begin{abstract}
On majore les noyaux de Green et de la chaleur au dehors de la diagonale pour des opérateurs
de type laplacien sur les espaces symétriques.
\end{abstract}
\maketitle

\section{Introduction}
On considère ici un espace symétrique $X=G/K$ de type non compact. À une représentation unitaire ($\rho,
V$) de dimension finie de $K$, on associe le fibré vectoriel $G\times_{K} V$ au dessus de $X$,
dont l'espace des sections lisses s'identifie à
$$C^\infty(E)\simeq\left\{\sigma\in C^\infty(G,V),\ \forall g\in G,\ \forall k\in K\,:\, f(gk)=\rho\left(k^{-1}\right)f(g)\,\right\}$$
L'objet de cet article est un opérateur de type laplacien $G-$invariant agissant sur les sections de $E$ 
$$L=\nabla^*\nabla+V$$ où
$\nabla$ est une connexion hermitienne $G$ invariante sur $E$ et $R$ un section $G-$invariante du fibré
des endomorphismes hermitiens de $E$. Nous donnons quelques estimations de la résolvante et du noyau de
la chaleur associé à $L$. Notre premier résultat est le suivant :
\begin{theo}\label{Green}Notons $\alpha_0$ le bas du spectre de l'opérateur $L$, $e=\Id.K\in X$ et  $G_s(x,y)$ le noyau de Schwarz de
la résolvante $\left(L-\alpha_0+s^2\right)^{-1}$ où $s$ est un nombre complexe de partie réelle strictement
positive. Il y a alors une constante $C$ telle que pour tout $x\in X$ tel que $d(x,o)\ge 2$ alors
$$|G_s(x,o)|\le C e^{-\rho(x^+)-\Re(s) d(x,o)}\ ,$$
où on a noté $x^+$ la composante suivant $\bar \Ag^+$ de $x$ dans la décomposition de Cartan
$G=Ke^{\bar \Ag^+}K$ et $\rho\in \Ag^*_\C$ la demi somme des racines restreintes positives associées à
$(\Gg_\C,\Ag)$.
\end{theo} 

Notre résultat est sensiblement meilleur que celui obtenu récemment par N. Lohoué et S. Mehdi à propos du
laplacien de Hodge deRham ; dans \cite{LM}, en utilisant un théorème de Paley-Wiener de P. Delorme
(\cite{D}) et la théorie des représentations de $G$, ils obtiennent pour tout $\epsilon\in ] 0,1 [$
l'existence d'une constante $C_\epsilon$ telle que pour tout $x\in X$ tel que $d(x,o)\ge 1$,
$$|G_s(x,o)|\le C_\epsilon \Phi_0(x)\, e^{-(1-\epsilon)\Re(s) d(x,o)}\ ,$$
où $\Phi_0$ est la fonction sphérique élémentaire de Harish-Handra de $G$.  On sait qu'il y a une
constante telle que $\Phi_0(x)\ge C\, e^{-\rho(x^+})$, en fait la fonction 
$\Phi_0(x) e^{\rho(x^+)}$ croit polynomialement sur $\bar \Ag^+$ (\cite{A}).

Cependant notre estimation n'est pas, en général, optimale. Par exemple pour les fonctions, on sait grâce
au travail de J-P. Anker et L. Ji (\cite[theorem 4.22 i)]{AJ}) que pour $s>0$, l'on a une estimation de la forme 
$$C^{-1} d(x,o)^{-\beta} \Phi_0(x) e^{-\rho(x^+)}\le G_s(x,o)\le 
C\, d(x,o)^{-\beta} \Phi_0(x) e^{-sd(x,o)}$$
où si on note $\Sigma^{++}$ les racines positives indivisibles et $l$ le rang de $X$ alors
$\beta=|\Sigma^{++}|+\frac{l-1}{2}$.
En fait, on a l'estimation
$d(x,o)^{-\beta} \Phi_0(x)\le Cd(x,o)^{-\frac{l-1}{2}}  e^{-\rho(x^+)}.$
Sur les fonctions, notre estimation est donc optimale en rang $1$ et en rang supérieur, 
elle est optimale à un facteur
polynomiale près  ; notons également que grâce à \cite[theorem 3.6]{CP}, notre résultat est optimal
 pour le laplacien de Hodge-deRham en rang $1$.

La preuve de notre estimation n'utilise que deux ingrédients à savoir un argument classique de
propagation à vitesse finie et une estimation du volume de $KB(x,1)\subset X$. Nous avons également obtenu
une estimation du noyau de Schwarz de l'opérateur de la chaleur $e^{-tL}$ par la même méthode. Pour
énoncer ce résultat, on rappelle quelques notations sur la structure algébrique de $X$. On note
$\Kg\subset\Gg$ les algèbres de Lie de $K$ et $G$ et $$\Gg=\Kg\oplus \Pg$$
 la décomposition en espaces propres de l'involution de Cartan. Si $\Ag\subset \Pg$ est une sous-algèbre
 abélienne maximale et $\Sigma\subset \Ag_\C^*$ le système restreint des racines de $(\Gg,\Ag)$. On fixe
 alors $\Ag^+\subset \Ag$ une chambre de Weil et on note $\Sigma^+\subset \Sigma$ le système des racines restreintes
 positives associés. Le rang de l'espace symétrique $X$ est $l=\dim \Ag$ ; la dimension de l'espace
 symétrique $X$ est notée $d$. L'espace $\Pg$ se décompose en
 $$\Pg=\Ag\oplus \Ng=\Ag\oplus\bigoplus_{\alpha\in \Sigma^+} \Ng_ \alpha.$$
 où pour  
 $a\in \Ag$ et $n\in \Ng_ \alpha$, nous avons $\ad(a)n=\alpha(a) n$.
 On note aussi $m_\alpha=\dim\Ng_ \alpha $ et donc $\rho=\frac12 \sum_{\alpha\in \Sigma^+}
 m_\alpha\alpha\in \Ag_\C^*$. Dans la décomposition de Cartan de $G=Ke^{\bar \Ag^+}K$, si $x=gK\in X=G/K$
, on note $x^+$  un élément de $\bar \Ag^+$ telle que $gK\in Ke^{x^+}K.$
Notre estimation est alors la suivante :
\begin{theo}\label{heat}Notons $h_t(x,y)$ le noyau de Schwarz de l'opérateur de la chaleur $e^{-tL}$. 
Il existe une constante $C$ telle que pour tout $x\in X$ tel que $d(x,e)\ge 2$ alors
$$|h_t(x,e)|\le C\, e^{-\alpha_0 t-\rho(x^+)-\frac{d(x,o)^2}{4t} }\phi_t(x)$$ où
$$\phi_t(x)=\left\{ \begin{array}{ll} \frac{ \sqrt{t}}{d(x,o)+\sqrt{t}}& \mbox{ si } d(x,o)\le t\\
\frac{d(x,o)^{\frac{d+l}{2}-1}}{t^{\frac{d+l-1}{2}}}\, \prod_{\alpha\in \Sigma^+}\left( 
\frac{1+\alpha(x^+)}{\frac{t}{d(x,0)}+\alpha(x^+)}
\right)^{m_\alpha/2}&\mbox{ si } d(x,o)\ge t\\
\end{array}\right.$$
\end{theo} 

Cette majoration n'est également pas optimale. On peut comparer notre estimation avec celle obtenue par 
N. Lohoué et S. Mehdi à propos du
laplacien de Hodge deRham  \cite{LM} ; ils obtiennent pour tout $\epsilon\in ]0,
1 [ $ des constantes $C_\epsilon$ et $A_\epsilon$ telles que  si $d(x,o)\ge A_\epsilon$ alors
$$|h_t(x,o)|\le C_\epsilon \Phi_0(x) e^{-\alpha_0 t} e^{-(1-\epsilon) \frac{d(x,o)^2}{4t} }
t^{-\epsilon\gamma}.$$
Notre estimation est donc meilleure lorsque $d(x,o)$ tend vers l'infini mais bien plus mauvaise lorsque
$t$ tend vers $+\infty$. Dans le cas de l'espace hyperbolique réel et du laplacien scalaire, 
on peut vérifier avec l'estimation de E. Davies et N. Mandouvalos \cite{DM} que notre estimée est optimale dans le régime où 
$d(x,o)\ge \max\{2,t\}$.

\section{Une estimée de volume}
Nous démontrons ici le résultat suivant :
\begin{prop}\label{volume} Il y a des constantes strictement positives $c,C$ telle que pour tout $\epsilon\in [0,1[$ et
$x\in X$ alors
$$ c e^{2\rho(x^+)}\, \epsilon^l\,  \prod_{\alpha\in \Sigma^+}\left( 
\frac{\epsilon+\alpha(x^+)}{1+\alpha(x^+)}
\right)^{m_\alpha}\le \vol KB(x,\epsilon)\,\le C e^{2\rho(x^+)}\, \epsilon^l\, \prod_{\alpha\in \Sigma^+}\left( 
\frac{\epsilon+\alpha(x^+)}{1+\alpha(x^+)}
\right)^{m_\alpha}.$$
\end{prop}

\proof Grâce à \cite[lemme 2.1.2]{AJ}, nous savons que
$$KB(x,\epsilon)\simeq K\exp(B(x^+,\epsilon)\cap\bar\Ag^+)$$ 
dans la décomposition de Cartan $X=Ke^{\bar\Ag^+}$.
Ainsi si $J(h)dkdh$ est la forme volume de $X$ dans les coordonnées $(k,h)\mapsto ke^hK$
nous avons :
$$\vol KB(x,\epsilon)=\int_ {B(x^+,\epsilon)\cap\Ag^+} J(h)dh.$$
Cependant  nous avons 
$$J(h)=C \prod_{\alpha\in \Sigma^+} \sinh^{m_\alpha}(\alpha(h))\approx e^{2\rho(h)} 
\prod_{\alpha\in \Sigma^+}\left( 
\frac{\alpha(h)}{1+\alpha(h)}\right)^{m_\alpha}.$$

Grâce à la preuve de \cite[lemme 2.1.6 i)]{AJ}, on en déduit que pour  $\epsilon\in
]0,1[$ et $h\in B(x^+,\epsilon)\cap\Ag^+$, on a
$$\rho(x^+)-|\rho|\le \rho(h)\le \rho(x^+)+|\rho|$$
et pour $\alpha\in \Sigma^+$, 
$$|\alpha(h-x^+)|\le \epsilon/\sqrt{2}$$
\begin{equation*}
\begin{split}
&\left(1-\frac{1}{\sqrt{2}}\right) (1+\alpha(x^+))\le 1+\alpha(h)\le 2 (1+\alpha(x^+))\\
&\alpha(h)\le \alpha(x^+)+\epsilon
\end{split}\end{equation*}
On en déduit aisément la majoration annoncée.

Pour la minoration, on considère
$\Sigma^{+++}=\{\alpha_1,...,\alpha_l\}$ un système
de racines réduites qui est une base de $\Ag_\C^*$
 et $E_1,...,E_l$ la base de $\Ag$ duale à $\Sigma^{+++}$,  si
  $\sum_{\alpha\in \Sigma^+}\alpha = \sum_{i=1}^l x_i \alpha_i$, on pose alors 
 $$v=\sum_i  E_i.$$ Ainsi pour $\alpha\in \Sigma^+$, on a $\alpha(v)\ge 1$.
 On a bien évidemment 
 $$B\left(x^++\frac{\epsilon}{10+10|v|}v,\frac{\epsilon}{20+20|v|} \right)\subset B(x^+,\epsilon)$$
 or sur $B\left(x^+ +\frac{\epsilon}{10+10|v|}v,\frac{\epsilon}{20+20|v|}\right) $, 
 on a pour  $\alpha\in \Sigma+$
 $$\alpha\ge \alpha(x^+)+\frac{\epsilon}{10+10|v|}\alpha(v)-\frac{\epsilon}{20+20|v|}\ge\alpha(x^+)+\frac{\epsilon}{20+20|v|} .$$
  On obtient ainsi facilement une minoration du volume de
  $KB\left(x^++\frac{\epsilon}{10+10|v|}v,\frac{\epsilon}{20+20|v|} \right)$ et donc du volume de
  $KB(x^+,\epsilon)$.
\endproof
\section{Estimation du noyau de Green} Ici, on étudie le noyau de Schwarz de l'opérateur
$(L-\alpha_0+s^2)^{-1}$ au dehors de la diagonale où $s$ est un nombre complexe de partie réelle strictement
positive. On commence par une estimée classique induite par la
propriété de propagation à vitesse finie de l'opérateur $\cos \left(t\sqrt{L-\alpha_0}\,\right)$ (cf. par exemple
 \cite[appendice D]{MM}).
 \begin{lem} Soit $\sigma\in L^2(B(o,1),E)$ et $u:=(L-\alpha_0+s^2)^{-1}\sigma$ alors pour $A:=KB(x,1)$
 avec $x\in X$ vérifiant $d(x,o)\ge 2$, on a 
 $$\|u\|_{L^2(A)}\le \frac{1}{(\Re s)^2} e^{-\Re(s)(d(x,o)-2)}\, \|\sigma\|_{L^2}.$$
 \end{lem}
 \proof  En effet, on a 
 $$u=\int_0^\infty \frac{e^{s\xi}}{s} \cos\left(\xi\sqrt{L-\alpha_0}\,\right)\sigma\, d\xi .$$
 Les hypothèses faites sur $x$ et $\sigma$ et la propriété de propagation à vitesse finie implique que
dès que $0\le \xi\le d(x,o)-2$ , on a 
 $\ \|\cos\left(\xi\sqrt{L-\alpha_0}\,\right)\sigma\|_{L^2(A)}=0$ .
 D'où
\begin{equation*}
\begin{split}
\|u\|_{L^2(A)}&\le\int_{d(x,o)-2}^\infty \frac{e^{-\Re(s)\xi}}{|s|}
\left\|\cos\left(\xi\sqrt{L-\alpha_0}\,\right)\sigma\right\|_{L^2(A)}\ d\xi\\
&\le \int_{d(x,o)-2}^\infty \frac{e^{-\Re(s)\xi}}{|s|}
\left\|\cos\left(\xi\sqrt{L-\alpha_0}\,\right)\sigma\right\|_{L^2(X)}\ d\xi\\
&\le \int_{d(x,o)-2}^\infty \frac{e^{-\Re(s)\xi}}{|s|}
\|\sigma\|_{L^2(X)}\ d\xi\\
&\le\frac{1}{(\Re s)^2} e^{-\Re(s)(d(x,o)-2)} \|\sigma\|_{L^2}.
\end{split}
\end{equation*}
  \endproof
  \begin{lem}\label{bonneboule}Il y a une constante $C$ qui ne dépend que que $X$ et $s$ telle qu'il y a un $k\in K$ tel que
  $$\|u\|_{L^2(B(kx,\frac14))}\le C e^{-\Re(s)d(x,o)-\rho(x^+)} \|\sigma\|_{L^2}.$$
  \end{lem} 
  \proof En effet
 \begin{equation*}\begin{split}
 \int_{KB(x,\frac14)} \left(\int_{B(y,\frac12)} |u|^2(z)dz\right)dy&=
 \int_{KB(x,1)} \vol\left(B\left(z,1/2\right)\cap KB\left(x,1/4\right)\right)\, |u|^2(z)dz\\
 &\le C\int_{KB(x,1)}|u|^2(z)dz.
 \end{split}
\end{equation*}
On en déduit donc l'existence d'un $y\in KB\left(x,\frac14\right)$ tel que
$$\left(\vol KB\left(x,\frac14\right)\right)^{1/2} \|u\|_{L^2\left(B(y,1/2)\right)}\le C e^{-\Re(s)d(x,o)} \|\sigma\|_{L^2}.$$
Il y a donc un $k\in K$ telle que $d(y,kx)\le \frac14$ et d'où
$B(kx,\frac14)\subset B(y,\frac12)$ et grâce à l'estimée (\ref{volume}) :  
$\vol(KB(x,\frac14))^{1/2}\approx e^{\rho(x^+)}$, on obtient bien le résultat annoncé.
  \endproof
  
On utilise alors l'estimée elliptique standard suivante :
\begin{prop} Soit $\lambda \in \C$, il y a une constante $C$ qui dépend de $X,\lambda,L$ telle que si $r\in
]0,1]$ et si $v\in L^2(B(x,r),E)$ vérifie $Lv=\lambda v$ alors
$$|v(x)|\le \frac{C}{r^{d/2}} \|v\|_{L^2(B(x,r))}.$$
\end{prop} 
On en déduit donc l'estimation :
$$|u(kx)|\le Ce^{-\Re(s)d(x,o)-\rho(x^+)} \|\sigma\|_{L^2}.$$
Ceci implique que
$$\left(\int_{B(o,1)} |G_s(kx,y)|^2 dy \right) ^{1/2}\le Ce^{-\Re(s)d(x,o)-\rho(x^+)} .$$
Les mêmes estimées elliptiques fournissent alors la majoration suivante :
$$\mbox{pour }d(x,o)\ge 2\ :\ |G_s(kx,o)|=|G_s(x,o)|\le Ce^{-\Re(s)d(x,o)-\rho(x^+)}.$$
Ce qui démontre le théorème (\ref{Green}).
\section{Estimation du noyau de la chaleur}
On étudie maintenant de la même façon le noyau de Schwarz de l'opérateur de la chaleur $e^{-tL}$ au
dehors de la diagonale. Nous commençons par le même type d'estimations :
\begin{lem}
Soit $\sigma\in L^2(B(o,\epsilon),E)$ et $f_t:=e^{-tL}\sigma$ alors pour $A:=KB(x,\epsilon)$
 avec $x\in X$ vérifiant $d(x,o)\ge 2\epsilon$, on a 
 $$\|f_t\|_{L^2(A)}\le \frac{e^{-\alpha_0 t}}{  \sqrt{\pi t}} \int_{d(x,o)-2\epsilon}^\infty
 e^{-\frac{\xi^2}{4t}} d\xi\,\|\sigma\|_{L^2}.$$
 \end{lem}
 Cette estimation se montre de la même façon que l'estimée en partant de la formule :
 $$f_t=
 \frac{e^{-\alpha_0 t}}{  \sqrt{\pi t}}\int_0^\infty e^{-\frac{\xi^2}{4t}}
 \cos\left(\xi\sqrt{L-\alpha_0}\,\right)\sigma\, d\xi.$$
 En reprenant la preuve du lemme (\ref{bonneboule}),
  on obtient que si de plus $\epsilon\in ]0,1/2]$ alors on trouve un $k\in
 K$ tel que
 
$$\|f_t\|^2_{L^2(B(kx,\frac\epsilon4))}\le C \frac{\epsilon^d}{\vol KB(x,\epsilon)} \|f_t\|^2_{L^2(A)}.$$
Maintenant, on utilise l'estimation parabolique suivante :
\begin{prop} Il y a une constante $C$ (qui ne dépend que de $X,L$) telle que si $r\in ]0,1]$ et si 
$v\in L^2([t-r^2,t]\times B(x,r), E)$ est une solution de l'équation :
$$\frac{\partial}{\partial t} v+Lv=0$$ alors
$$|v(t,x)|^2\le \frac{C}{r^{d+2}}\int_{t-r^2}^t\left(\int_{B(x,r)} |v(\tau,y)|^2 dy \right)d\tau.$$
\end{prop}

Il y a donc une constante $C$ tel que
pour $\epsilon\in ]0,1/2]$ et $t\ge 2\epsilon^2$ alors
$$|f_t(kx)|^2\le \frac{C}{\epsilon^{d+2}} \int_{[t-\epsilon^2,t]\times B(kx,\frac\epsilon4)} |f_\tau(y)|^2
d\tau dy.$$
Avec l'estimation 
$$\int_A^\infty e^{-\frac{\xi^2}{4t}} d\xi=\sqrt{t} \int_{A^2/4t}^\infty e^{-v} \frac{dv}{\sqrt{v}}\le 
\frac{C\sqrt{t}}{\frac{A}{\sqrt{t}}+1} e^{-A^2/4t}\ ,$$
on obtient, pour $t\ge 2\epsilon^2$ et $d(x,o)\ge 2$, l'estimée suivante 
$$|f_t(kx)|\le 
 C\,(\vol KB(x,\epsilon))^{-1/2}\frac{ \sqrt{t}}{d(x,o)+\sqrt{t}}\,
 e^{-\alpha_0 t-\frac{ (d(x,o)-2\epsilon)^2}{4t}}\; \|\sigma\|_{L^2}\ .$$
 
 On en déduit 
 que
 $$\left(\int_{B(o,\epsilon)} |h_t(kx,y)|^2 dy \right) ^{1/2}\le C\,(\vol KB(x,\epsilon))^{-1/2}\frac{ \sqrt{t}}{d(x,o)+\sqrt{t}} 
 \,e^{-\alpha_0 t-\frac{ (d(x,o)-2\epsilon)^2}{4t}}\ .$$
 
 Les mêmes estimées paraboliques donnent alors la majoration suivante : pour $\epsilon\in ]0,1/2]$, $t\ge
 3\epsilon^2$ et $d(x,o)\ge 2$ :
 $$|h_t(kx,o)|=|h_t(x,o)|\le C (\epsilon^{d}\vol KB(x,\epsilon))^{-1/2}\frac{ \sqrt{t}}{d(x,o)+\sqrt{t}} 
 \,e^{-\alpha_0 t-\frac{ (d(x,o)-2\epsilon)^2}{4t}}.$$
 Or nous avons
  $$\frac{(d(x,o)-2\epsilon)^2}{4t}=\frac{d(x,o)^2}{4t}-\frac{d(x,o)\epsilon}{t}+\frac{\epsilon^2}{t}.$$
  Donc lorsque $d(x,o)\le t$ on choisit $\epsilon=1/100$ et on obtient la majoration :
$$ |h_t(x,o)|\le C\,\frac{ \sqrt{t}}{d(x,o)+\sqrt{t}} \, 
 e^{-\alpha_0 t -\frac{d(x,o)^2}{4t}-\rho(x^+)}$$
Lorsque $d(x,o)\ge t$ on choisit $\epsilon=\frac{t}{100d(x,o)}$ et on obtient pour $d(x,o)\ge 2$
 \begin{equation*}
 \begin{split}
 |h_t(x,o)|&\le C \frac{\sqrt{t}}{d(x,o)} \left(\frac{d(x,o)}{t}\right)^{\frac{d+l}{2}}\prod_{\alpha\in
\Sigma^+}\left(\frac{1+\alpha(x^+)}{\frac{t}{100d(x,o)}+\alpha(x^+)}\right)^{m_\alpha/2}\ 
 e^{-\alpha_0 t-\frac{d(x,o)^2}{4t}-\rho(x^+)}\\
 &\le C \frac{d(x,o)^{\frac{d+l}{2}-1}}{t^{\frac{d+l-1}{2}}}\prod_{\alpha\in
\Sigma^+}\left(\frac{1+\alpha(x^+)}{\frac{t}{d(x,o)}+\alpha(x^+)}\right)^{m_\alpha/2}\ 
 e^{-\alpha_0 t-\frac{d(x,o)^2}{4t}-\rho(x^+)}\end{split}\end{equation*}
 Ce qui termine la démonstration du théorème (\ref{heat}).
 \section{Applications} Dans \cite{CP}, une estimation du prolongement
 analytique de la résolvante avait été obtenue  ; cependant les méthodes rudimentaires utilisées
  ici ne permettent pas d'obtenir un tel résultat. Néanmoins, nos estimées comme
   celle de N.Lohoué et S. Mehdi permettent une
 estimation inférieure du bas du spectre de l'opérateur $L$ sur des espaces localement symétriques
 $\Gamma\backslash G/K$ où $\Gamma\subset G$ est un sous-groupe discret sans torsion.
 
\noindent{\bf Définition :} (cf. \cite[theorem 2.7]{CP})
Soit $\Gamma\subset G$ est un sous-groupe discret sans torsion, on note $\tilde \delta(\Gamma)$
l'exposant critique modifié de $\Gamma$ qui est l'exposant critique de la série de Poincaré:
$$\sum_{\gamma\in \Gamma} e^{\rho(\gamma^+)-sd(\gamma(o),o)}\ .$$  
Notons $G^0_s(x,y)$ le noyau de Green de l'opérateur $(\Delta-|\rho|^2+s^2)^{-1}$ agissant sur les fonctions. 
Grâce à notre estimation et à l'estimation inférieure de $G^0$ de J-P. Anker et L.Ji 
(\cite[Theorem 4.2.2]{AJ}), on sait que pour tout $s>0$ et 
$\eta\in ]0,s[$, il y a une constante $C_{s,\eta}$ tel que pour tout $x,y\in X$, 
$$|G_s(x,y)|\le C_{s,\eta} G^0_{s-\eta}(x,y)$$
Le même raisonnement que celui utilisé pour démontrer \cite[theorem 2.7]{CP} montre que :
\begin{thm}\begin{enumeroman}
\item Si $\tilde \delta(\Gamma)>0$ alors le bas du spectre de $L$ sur  $\Gamma\backslash G/K$ est
minoré par $\alpha_0-(\tilde\delta(\Gamma) )^2$.

\item Si $\tilde\delta(\Gamma)\le 0$ alors le bas du spectre de $L$ sur  $\Gamma\backslash G/K$ est
minoré par $\alpha_0$. 
\item Si  $\tilde\delta(\Gamma)\le 0$ et si le rayon d'injectivité de $\Gamma\backslash G/K$ est non-majoré i.e.
$$\sup_{x\in X}\inf_{\gamma\in \Gamma}d(x,\gamma(x))=\infty$$ alors le bas du spectre de $L$ est 
$\alpha_0$. 
\end{enumeroman}
\end{thm}
\begin{rem}\begin{enumeroman}
\item Cet exposant critique modifié se compare aisément à l'exposant critique de $\Gamma$, à savoir
à $\delta(\Gamma)$ l'exposant critique de la série 
$$\sum_{\gamma\in \Gamma} e^{-sd(\gamma(o),o)}.$$
Si on note $\rho_{\min}=\inf_{H\in \Ag_+} \rho(h)/|h|$ alors
$$\rho_{\min}+\tilde\delta(\Gamma)\le \delta(\Gamma)\le |\rho|+\tilde\delta(\Gamma)$$
ce qui permet de ré-obtenir le résultat de N. Lohoué et S. Mehdi \cite[theorem 6.1]{LM}.

\item Lorsque $\tilde\delta(\Gamma)< \sqrt{\alpha_0}$, ce résultat implique que le noyau $L^2$ de $L$ sur
$\Gamma\backslash G/K$ est trivial. Il est cependant difficile d'obtenir des calculs explicites de 
$\alpha_0$. Concernant le laplacien de Hodge-deRham sur les formes différentielles des calculs explicites
sont fait par H. Donnelly, E. Pedon  en rang $1$ (\cite{Do,P1,P2}) et par  N. Lohoué et S. Medhi pour certains
espaces hermitiens \cite[Appendix A]{LM}.
\end{enumeroman}\end{rem}

\end{document}